\documentclass[11pt]{article}
\usepackage{amssymb,amsfonts,amsmath,amsthm}
\usepackage{epsfig}
\newtheorem{theorem}{Theorem}[section]
\newtheorem{thm}[theorem]{Theorem}

\newtheorem{lem}[theorem]{Lemma}

\newtheorem{cor}[theorem]{Corollary}

\newtheorem{conj}[theorem]{Conjecture}

\makeatletter \@addtoreset{equation}{section}

\def\CT{\mathop{\mathrm{CT}}}

\begin{document}

\title{On Kadell's two Conjectures for the $q$-Dyson Product}
\author{
  {\vspace{0.2cm}}
  {Yue Zhou}\\
 {\small School of Mathematical Science and Computing Technology}\\
{\small Central South University, Changsha 410075, P.R. China}\\
  {nkzhouyue@gmail.com}\\
 }

\date{September 7, 2010}
\maketitle

\begin{abstract}
By extending Lv-Xin-Zhou's first layer formulas of the $q$-Dyson
product, we prove Kadell's conjecture for the Dyson product and show
the error of his $q$-analogous conjecture. With the extended
formulas we establish a $q$-analog of Kadell's conjecture for the
Dyson product.
\end{abstract}



\section{Introduction}

In 1962, Freeman Dyson \cite{dyson1962} conjectured the following
constant term identity.
\begin{theorem}[Dyson's Conjecture]\label{t-dyson}
For nonnegative integers $a_0,a_1,\ldots ,a_n$,
\begin{equation*}
\CT_{\mathbf{x}} \prod_{0\leqslant i\ne j \leqslant n}
\left(1-\frac{x_i}{x_j}\right)^{\!\!a_i} =
 \frac{a!}{a_0!\, a_1!\, \cdots a_n!},
\end{equation*}
where $a:=a_0+a_1+\cdots+a_n$ and $\CT_{\mathbf{x}}f(\mathbf{x})$
means to take constant term in the $x$'s of the series
$f(\mathbf{x})$.
\end{theorem}

The conjecture was quickly proved independently by
Gunson~\cite{gunson} and by Wilson \cite{wilson}. An elegant
recursive proof was published by Good \cite{good1}, and a
combinatorial proof was given by Zeilberger \cite{zeil1982}. In
1975, George Andrews \cite{andrews1975} came up with a $q$-analog of
the Dyson conjecture.
\begin{thm}\label{thm-dyson}\emph{(Zeilberger-Bressoud)}.
For nonnegative integers $a_0,a_1,\dots,a_n$,
\begin{align*}
\CT_{\mathbf{x}}\,\prod_{0\leqslant i<j\leqslant n}
\left(\frac{x_i}{x_j}\right)_{\!\!a_i}
\left(\frac{x_j}{x_i}q\right)_{\!\!a_j}
=\frac{(q)_{a}}{(q)_{a_0}(q)_{a_1}\cdots(q)_{a_n}},
\end{align*}
where $(z)_m:=(1-z)(1-zq)\cdots(1-zq^{m-1})$.
\end{thm}
The Laurent polynomials in the above two theorems are respectively called
the \emph{Dyson product} and the \emph{$q$-Dyson product} and respectively denoted
by $D_n(\mathbf{x},\mathbf{a})$ and $D_n(\mathbf{x},\mathbf{a},q)$,
where $\mathbf{x}:=(x_0,\ldots,x_n)$ and
$\mathbf{a}:=(a_0,\ldots,a_n)$.

The Zeilberger-Bressoud $q$-Dyson Theorem was first proved,
combinatorially, by Zeilberger and Bressoud \cite{zeil-bres1985} in
1985. Recently, Gessel and Xin \cite{gess-xin2006} gave a very
different proof by using the properties of the formal Laurent series and of the
polynomials. The coefficients of the Dyson and the $q$-Dyson product
were researched in \cite{breg, kadell1998, Xin, Xin2, sills2006,
sills-zeilberger, stembridge1987}. In the equal parameter case, the
identity reduces to Macdonald's constant term conjecture
\cite{macdonald} for root systems of type $A$. In 1988 Stembridge \cite{stembridge-qdyson}
gave the first layer formulas of the $q$-Dyson product in the equal
parameter case.

Let $I=\{i_1,\ldots,i_{m}\}$ be a proper subset of
$\{0,1,\ldots,n\}$ and $J=\{j_1,\ldots,j_m\}$ be a multi-subset of
$\{0,1,\ldots,n\}\setminus I$, where $0\leqslant i_1<\cdots <i_m\leqslant n$
and $0\leqslant j_1\leqslant \cdots \leqslant j_{m}\leqslant n$.

Our first objective in this paper is to prove the following
conjecture of Kadell \cite{kadell1998}.

\begin{conj}\label{conj-Kadell}
For nonnegative integers $a_0,a_1,\ldots,a_n$ we have
\begin{align}\label{e-conj} \Big(1+a-\sum_{k\in I}a_k\Big)\CT_{\mathbf{x}} \prod_{k=1}^m\Big(1-\frac{x_{j_k}}{x_{i_k}}\Big)
\prod_{0\leqslant i\ne j \leqslant n}\left(1-\frac{x_i}{x_j}\right)^{\!\!a_i}
=\Big(1+a\Big)\frac{a!}{a_0!a_1!\cdots a_n!}.
\end{align}
\end{conj}
In the same paper, Kadell also gave a $q$-analogous conjecture, we
restate it as follows.
\begin{conj}\label{conj-Kadell2}
Let $P=\{(i_k,j_k)\mid i_{k}\in I,j_k\in J,k=1,2,\ldots,m\}$. Then for nonnegative integers $a_0,a_1,\ldots,a_n$ we have
\begin{align}\label{e-conj2} \Big(1-&q^{1+a-\sum_{k\in I}a_k}\Big)\CT_{\mathbf{x}} \prod_{0\leqslant s<t\leqslant
n}\left(\frac{x_s}{x_t}\right)_{\!\!a_i+\chi((t,s)\in P)}
\left(\frac{x_t}{x_s}q\right)_{\!\!a_j+\chi((s,t)\in P)} \nonumber \\
&=\Big(1-q^{1+a}\Big)\frac{(q)_{a}}{(q)_{a_0}(q)_{a_1}\cdots
(q)_{a_n}},
\end{align}
where the expression $\chi(S)$ is $1$ if the statement $S$ is true,
and $0$ otherwise.
\end{conj}


In trying to prove Conjecture \ref{conj-Kadell2}, we find that the
conjectured formula is incorrect. One way to modify the conjecture
is to evaluate the left-hand side of \eqref{e-conj2}. This can be
done by writing it as a linear combination of some first layer
coefficients of the $q$-Dyson product, and then applying the
formulas of \cite{Xin}. Unfortunately, we are not able to derive a
nice formula.

Our second objective is to contribute a $q$-analogous formula of
\eqref{e-conj}, which is motivated by the proof of \eqref{e-conj},
and is stated in Theorem \ref{t-m}.

This paper is organized as follows. In Section 2 we reformulate the
main result in \cite{Xin} and give an extended form of it. In
Section 3 we prove Conjecture \ref{conj-Kadell} and give an example
to show the error of Conjecture \ref{conj-Kadell2}.
In Section 4 based on Conjecture \ref{conj-Kadell2} we give our main theorem.

\section{Basic results}

Let $T=\{t_1,\ldots,t_d\}$ be a $d$-element subset of $I$ with
$t_1<\cdots<t_d$. Define
\begin{equation}\label{e-def-w}w_{i}=\left\{ \begin{aligned}
         &a_i,\quad &for& \quad i\not\in T; \\
                  &0, \quad &for& \quad i\in T.
                          \end{aligned} \right.\end{equation}
Let $S$ be a set and $k$ be an element in $\{0,1,\ldots,n\}$. Define
$N(k,S)$ as the number of the elements in $S$ which are not larger
than $k$, i.e.,
\begin{align}\label{def-N}
N(k,S)=\big|\{i\leqslant k\mid i\in S\}\big|.
\end{align}
In particular, $N(k,\varnothing)=0$.

The \emph{first layer formulas} of the $q$-Dyson product can be
restated as follows.
\begin{thm}\emph{\cite{Xin}}\label{t-main-thm}
Let $I,J$ be defined as in Conjecture \ref{conj-Kadell}. Then for
nonnegative integers $a_0,a_1,\ldots,a_n$ and fixed $i_1=0$ we have
{\small
\begin{align}\label{e-main}
\CT_{\mathbf{x}}\frac{x_{j_1}x_{j_2}\cdots
x_{j_{m}}}{x_{i_1}x_{i_2}\cdots
x_{i_m}}D_n(\mathbf{x},\mathbf{a},q) 
=\frac{(q)_{a}}{(q)_{a_0}\cdots(q)_{a_n}}\sum_{ \varnothing\neq
T\subseteq I}(-1)^dq^{L(T\mid I)} \frac{1-q^{\sum_{k\in
T}a_k}}{1-q^{1+a-\sum_{k\in T}a_k}},
\end{align}}
where
\begin{align}\label{e-LT}
L(T\mid I)=\sum_{k=0}^n\big[N(k,I)-N(k,J)\big]w_k.
\end{align}
\end{thm}
We need the explicit formula for the case $i_1\ne 0$ for our
calculation. As stated in \cite{Xin}, the formula for this case can
be derived using an action $\pi$ on Laurent polynomials:
\begin{align*}
\pi \big(F(x_0,x_1,\ldots,x_n)\big)=F(x_1,x_2,\dots,x_n,x_0/q).
\end{align*}

By iterating, if $F(x_0,x_1,x_2,\dots,x_n)$ is homogeneous of degree
$0$, then
$$\pi^{n+1}\big(F(x_0,x_1,\ldots,x_n)\big)=F(x_0/q,x_1/q,x_2/q,\ldots,x_n/q)=F(x_0,x_1,x_2,\ldots,x_n),  $$
so that in particular $\pi$ is a cyclic action on
$D_n(\mathbf{x},\mathbf{a},q)$. We use the following lemma to derive
an extended form of Theorem \ref{t-main-thm}.

\begin{lem}\emph{\cite{Xin}}\label{l-pi}
Let $L(\mathbf{x})$ be a Laurent polynomial
in the $x$'s. Then
\begin{align}\label{e-cyclic}
\CT_\mathbf{x} L(\mathbf{x})\,D_n(\mathbf{x},\mathbf{a},q)=
\CT_\mathbf{x} \, \pi
\big(L(\mathbf{x})\big)D_n\big(\mathbf{x},(a_n,a_0,\ldots,a_{n-1}),q\big).
\end{align}
By iterating \eqref{e-cyclic} and renaming the parameters,
evaluating $\CT_\mathbf{x}
L(\mathbf{x})\,D_n(\mathbf{x},\mathbf{a},q) $ is equivalent to
evaluating $\CT_\mathbf{x}
\pi^k(L(\mathbf{x}))\,D_n(\mathbf{x},\mathbf{a},q)$ for any integer
$k$.
\end{lem}

Assume for some $t$ we have $j_t<i_1$ and $j_{t+1}>i_1$. Let
$J^-=\{j_1,\ldots,j_t\}$ and $J^+=\{j_{t+1},\ldots,j_{m}\}$.
\begin{thm}\label{t-general}
For nonnegative integers $a_0,a_1,\ldots,a_n$ we have
{\small \begin{align}\label{e-main2}
\CT_{\mathbf{x}}\frac{x_{j_1}x_{j_2}\cdots
x_{j_{m}}}{x_{i_1}x_{i_2}\cdots
x_{i_m}}D_n(\mathbf{x},\mathbf{a},q) 
=\frac{(q)_{a}}{(q)_{a_0}\cdots(q)_{a_n}}\sum_{ \varnothing\neq
T\subseteq I}(-1)^dq^{L^*(T\mid I)} \frac{1-q^{\sum_{k\in
T}a_k}}{1-q^{1+a-\sum_{k\in T}a_k}},
\end{align}}
where {\small \begin{align}\label{e-LT2} L^{*}(T\mid
I)=t+\sum_{k=i_1}^n\big[N(k,I)-N(k,J^+)\big]w_k+\sum_{k=0}^{i_1-1}\big[t-N(k,J^-)\big]a_k.
\end{align}}
\end{thm}
We remark that there is not the restriction $i_1=0$ in the above theorem.
The idea to prove this theorem is by iterating
Lemma \ref{l-pi} to transform the random $i_1$ in \eqref{e-main2} to
zero and then applying Theorem \ref{t-main-thm}. But in the proof there are
many tedious transformations of the parameters, so we put the proof
to the appendix for those who are interested in.

Note that $w_k$ only occurs in the first summation of \eqref{e-LT2},
so only the first summation of \eqref{e-LT2} depends on $T$.

Letting $q\rightarrow 1^{-}$ in Theorem \ref{t-general} we get
\begin{cor}\label{t-1}\emph{\cite{Xin}}
For nonnegative integers $a_0,\ldots,a_n$ we have
\begin{equation}\label{e-firstlayer}{\small
\CT_{\mathbf{x}}\frac{x_{j_1}\cdots x_{j_{m}}}{x_{i_1}\cdots
x_{i_{m}}} \prod_{0\leqslant i\ne j \leqslant n}
\left(1-\frac{x_i}{x_j}\right)^{\!\!a_i} =\frac{a!}{ a_0!\, \cdots
a_n!}\sum_{\varnothing\neq T\subseteq I}(-1)^{d}\frac{\sum_{k\in
T}a_k}{1+a-\sum_{k\in T}a_k}.}
\end{equation}
\end{cor}
This result also follows from \cite[Theorem 1.7]{Xin} by permuting the
variables. Note that the right-hand side of \eqref{e-firstlayer} is independent of $j$'s.

\section{Proof of Conjecture \ref{conj-Kadell}}

Now we are ready to prove Conjecture \ref{conj-Kadell}.

\begin{proof}[Proof of Conjecture {\bf \ref{conj-Kadell}}]
If $I=\varnothing$ then Conjecture \ref{conj-Kadell} reduces to the
Dyson Theorem, which is also the case when $m=0$ in Corollary
\ref{t-1}. So we assume that $I\neq \varnothing$.
Expanding the first product of \eqref{e-conj} gives {\small
\begin{align*} \CT_{\mathbf{x}} \prod_{i=1}^m\Big(1-\frac{x_{j_k}}{x_{i_k}}\Big)
\prod_{0\leqslant i\ne j \leqslant n}\left(1-\frac{x_i}{x_j}\right)^{\!\!a_i}
=\CT_{\mathbf{x}}\bigg[1+\sum_{l=1}^m(-1)^{l}\sum_{\varnothing\neq I_{l}\subseteq I}
\frac{x_{v_1}\cdots x_{v_l}}{x_{u_1}\cdots x_{u_l}}\bigg]
\prod_{0\leqslant i\ne j \leqslant n}\left(1-\frac{x_i}{x_j}\right)^{\!\!a_i},
\end{align*}}
where $I_{l}=\{u_1,\ldots,u_l\}$ ranges over all subsets of $I$
except the empty set and $\{v_1,\ldots,v_l\}$ is the corresponding
subset of $J$. Denote the left constant term in the above equation
by $LC$.
Applying Corollary \ref{t-1}, we get
\begin{align}\label{e-p1}
LC=\bigg[1+\sum_{l=1}^m(-1)^{l}\sum_{\varnothing\neq I_{l}\subseteq
I} \sum_{\varnothing\neq T\subseteq I_l}(-1)^{d}\frac{\sum_{k\in
T}a_k}{1+a-\sum_{k\in T}a_k}\bigg] \frac{a!}{ a_0!\, \cdots a_n!},
\end{align}
where $d=|T|$.
Changing the order of the summations, and observing that for any
fixed set $T$ there are $m-d\choose l-d$ such $I_l$ satisfying
$T\subseteq I_l\subseteq I$, we obtain
\begin{align}
LC&=\bigg[1+\sum_{\varnothing \neq T\subseteq
I}\sum_{l=d}^m(-1)^{l+d}{m-d\choose l-d}\frac{\sum_{k\in
T}a_k}{1+a-\sum_{k\in T}a_k}\bigg] \frac{a!}{ a_0!\, \cdots
a_n!} \nonumber\\
&= \bigg(1+\frac{\sum_{k\in I}a_k}{1+a-\sum_{k\in
I}a_k}\bigg)\frac{a!}{ a_0!\, \cdots a_n!},\label{e-p2}
\end{align}
where we used the easy fact that
for $d\neq m$
$$\sum_{l=d}^m(-1)^{l+d}{m-d \choose l-d}=\sum_{l=0}^{m-d}(-1)^l{m-d\choose l}=(1-x)^{m-d}\big|_{x=1}=0.$$
The conjecture then follows by multiplying
both sides of \eqref{e-p2} by $1+a-\sum_{k\in I}a_k$.
\end{proof}

For the $q$-case, Conjecture \ref{conj-Kadell2} does not hold even
for $m=1$. To see this take $n=2, I=\{0\}, J=\{1\}$ and $a_0=a_1=a_2=1$. For these values the left-hand side of \eqref{e-conj2} is
\begin{multline*}
(1-q^3)\CT_{\mathbf{x}}(1-\frac{x_0}{x_1})(1-q\frac{x_1}{x_0})(1-q^2\frac{x_1}{x_0})(1-\frac{x_0}{x_2})
(1-q\frac{x_2}{x_0})(1-\frac{x_1}{x_2})(1-q\frac{x_2}{x_1})\\
=(1-q^3)(1+2q+3q^2+2q^3).
\end{multline*}
While the right-hand side of \eqref{e-conj2} equals
$(1-q^4)(1+q)(1+q+q^2)$, which is not equal to the left-hand side.
%
%

\section{A $q$-analog of Kadell's conjecture}

\subsection{Motivation and presentation of the main theorem}
In this section we will construct a $q$-analog of Conjecture
\ref{conj-Kadell}. The new identity is motivated by the proof of
Conjecture \ref{conj-Kadell} in the last section, where massive
cancelations happen. We hope for similar cancelations in the
$q$-case.

Our first hope is to modify Conjecture \ref{conj-Kadell2} to obtain
a formula of the form:
{\small
\begin{align}\label{e-guess1}
\Big(1-&q^{1+a-\sum_{k\in I}a_k}\Big)\CT_{\mathbf{x}}
\prod_{k=1}^m\big(1-q^{L_{k}}\frac{x_{j_k}}{x_{i_k}}\big)D_{n}(\mathbf{x},\mathbf{a},q)
=\Big(1-q^{1+a}\Big)\frac{(q)_{a}}{(q)_{a_0}(q)_{a_1}\cdots
(q)_{a_n}},
\end{align}}
where $L_k$ is an integer depending on $i_k,j_k$ and $\mathbf{a}$.

It is intuitive to consider the $m=2$ case, so take $I=\{i_1,i_2\}$.
We need to choose appropriate $L_1$ and $L_2$ such that {\small
\begin{align}\label{e-guess2}
\Big(1-&q^{1+a-a_{i_1}-a_{i_2}}\Big)\CT_{\mathbf{x}}
\big(1-q^{L_1}\frac{x_{j_1}}{x_{i_1}}\big)\big(1-q^{L_2}\frac{x_{j_2}}{x_{i_2}}\big)D_{n}(\mathbf{x},\mathbf{a},q)
=\Big(1-q^{1+a}\Big)\frac{(q)_{a}}{(q)_{a_0}(q)_{a_1}\cdots
(q)_{a_n}}.
\end{align}}
By applying Theorem \ref{t-general}, the left-hand side of \eqref{e-guess2} becomes {\small
\begin{eqnarray}\label{e-guess3}
\Big(1-q^{1+a-a_{i_1}-a_{i_2}}\Big)\Big(1+q^{L_1+L^*(\{i_1\}\mid
\{i_1\})}\frac{1-q^{a_{i_1}}}{1-q^{1+a-a_{i_1}}}
+q^{L_2+L^*(\{i_2\}\mid \{i_2\})}\frac{1-q^{a_{i_2}}}{1-q^{1+a-a_{i_2}}}\nonumber \\
-q^{L_1+L_2+L^*(\{i_1\}\mid \{i_1,i_2\})}\frac{1-q^{a_{i_1}}}{1-q^{1+a-a_{i_1}}}
-q^{L_1+L_2+L^*(\{i_2\}\mid \{i_1,i_2\})}\frac{1-q^{a_{i_2}}}{1-q^{1+a-a_{i_2}}}\nonumber \\
+q^{L_1+L_2+L^*(\{i_1,i_2\}\mid \{i_1,i_2\})}\frac{1-q^{a_{i_1}+a_{i_2}}}{1-q^{1+a-a_{i_1}-a_{i_2}}}\Big)
\frac{(q)_{a}}{(q)_{a_0}(q)_{a_1}\cdots (q)_{a_n}}.
\end{eqnarray}}
It is natural to have the following requirements to get
\eqref{e-guess2}.
\begin{align}
q^{L_1+L^*(\{i_1\}\mid \{i_1\})}-q^{L_1+L_2+L^*(\{i_1\}\mid
\{i_1,i_2\})}=0, \nonumber \\
q^{L_2+L^*(\{i_2\}\mid \{i_2\})}-q^{L_1+L_2+L^*(\{i_2\}\mid
\{i_1,i_2\})}=0,\label{e-restrict}\\
q^{L_1+L_2+L^*(\{i_1,i_2\}\mid
\{i_1,i_2\})}=q^{1+a-a_{i_1}-a_{i_2}}. \nonumber
\end{align}
This is actually a linear system and has no solution, so our first
hope broke.

Looking closer at \eqref{e-restrict}, we see that the first two
equalities must be satisfied to have a nice formula. Agreeing with
this, for general $I$ with $|I|=m$ we will need $2^m-1$ restrictions
for massive cancelations as in the proof of Conjecture
\ref{conj-Kadell}. More precisely, by applying Theorem
\ref{t-general}, the left-hand side of \eqref{e-guess1} will be
written as
$$\Big(1-q^{1+a-\sum_{k\in I}a_k}\Big)\bigg(1+\sum_{T}B_T \frac{1-q^{\sum_{k\in T}a_k}}{1-q^{1+a-\sum_{k\in T}a_k}}\bigg)\frac{(q)_a}{(q)_{a_0}\cdots
(q)_{a_n}} ,$$ where $T$ ranges over all subsets of $I$ except the empty set. We need to
have $B_T=0$ for all $T$ except for $T=I$. This
is why using only $m$ unknowns dooms to fail.

We hope for some nice $A_T$ such that the constant term of
$$\sum_{T}A_T \frac{x_{v_1}\cdots x_{v_l}}{x_{u_1}\cdots x_{u_l}}
D_n(\mathbf{x},\mathbf{a},q)  $$ has the desired cancelations. We
are optimistical because from the view of linear algebra, such $A_T$
exists but is difficult to solve and might only be rational in $q$.
Amazingly, it turns out that in many situations, the $A_T$ may be
chosen to be $\pm q^{integer}$. Our formula for $A_T$ is inspired by
the proof of Conjecture \ref{conj-Kadell}. To present our result, we
need some notations.

Fix a subset $I=\{i_1,\ldots,i_m\}$ and a multi-subset
$J=\{j_1,\ldots,j_m\}$ of $\{0,1,\ldots,n\}$, where
$i_1<\cdots<i_m$, $j_1\leqslant \cdots \leqslant j_m$ and $I\cap
J=\varnothing$, $0\leqslant m\leqslant n$. Given an $l$-element subset
$I_l=\{u_1,\ldots,u_l\}$ of $I$, we say $J_l=\{v_1,\ldots,v_l\}$ is
the \emph{pairing set} of $I_l$ if $u_k=i_t\ (1\leqslant k\leqslant l)$ for
some $t$ implies that $v_k=j_t$. Write $I\setminus
I_l=\{i_{r_1},\ldots,i_{r_{m-l}}\}$, $r_1<\cdots<r_{m-l}$. We use
$A\mathop{\longrightarrow}\limits^{i} B$ to denote $B=A\cup \{i\}$,
and define a sequence of sets:
\begin{align}\label{def-II}I_l=\mathbb{I}_{m-l+1}\mathop{\longrightarrow}\limits^{i_{r_{m-l}}} \mathbb{I}_{m-l}
\mathop{\longrightarrow}\limits^{i_{r_{m-l-1}}} \mathbb{I}_{m-l-1}\mathop{\longrightarrow}\limits^{i_{r_{m-l-2}}}\cdots
\mathop{\longrightarrow}\limits^{i_{r_{1}}}\mathbb{I}_1=I.
\end{align}
For a set $S$ of integers, we denote by $\min S$ the smallest
element of $S$. Define $J_{k}^*(J_l)$ to be the set $\{j_s>\min
\mathbb{I}_k \mid j_s\in J_l\cup \{j_{r_k}\}\}$, we use $J_{k}^*$
as an abbreviation for $J_{k}^*(J_l)$.

Our $q$-analog of Conjecture \ref{conj-Kadell} can be stated as
follows.
\begin{thm}\emph{(Main Theorem)}\label{t-m}
For nonnegative integers $a_0,a_1,\ldots,a_n$ and $I,J$ as above, if
there is no $s,t,u$ such that $1\leqslant s<t<u\leqslant m$ and
$j_t<i_s<j_u<i_t$, then
 {\small
\begin{align}\label{e-last} \Big(1-&q^{1+a-\sum_{k\in I}a_k}\Big)\CT_{\mathbf{x}}\Bigg[\bigg(1+ \sum_{\varnothing \neq I_l\subseteq
I}(-1)^lq^{C(I_l)}\frac{x_{v_1}\cdots x_{v_l}}{x_{u_1}\cdots x_{u_l}}\bigg)
D_n(\mathbf{x},\mathbf{a},q)\Bigg] \nonumber \\
&=\Big(1-q^{1+a}\Big)\frac{(q)_{a}}{(q)_{a_0}(q)_{a_1}\cdots
(q)_{a_n}},
\end{align}}
where, with  $L^*(I_l\mid I_l)$ defined as in \eqref{e-LT2}, {\small
\begin{align}\label{e-CIl}
C(I_l)=1+a-\sum_{k\in
I_l}a_k+\sum_{k=1}^{m-l}\big[N(i_{r_k},I_l)-N(i_{r_k},J_{k}^*)\big]a_{i_{r_k}}-L^*(I_l\mid
I_l).
\end{align}}
\end{thm}
We remark that there is no
  analogous simple formula if the $u$'s and the $v$'s are not paired up,
  and that  
the sum $1+ \sum_{\varnothing \neq I_l\subseteq
I}(-1)^lq^{C(I_l)}\frac{x_{v_1}\cdots x_{v_l}}{x_{u_1}\cdots
x_{u_l}}$ in \eqref{e-last} does not factor.

\subsection{Factorization and cancelation lemma}

To prove the main theorem, we need some lemmas.

Let $U$ be a subset of $I_l$, $|U|=d$ and $I\setminus
U=\{i_{t_1},\ldots,i_{t_{m-d}}\}$, $t_1<\cdots<t_{m-d}$.
For fixed $I_l$, suppose that $\min I_l=i_v$.
By tedious calculation we can get the following lemma.
\begin{lem}\label{lem-subtract}
Let $U, C(I_l),L^*(U\mid I_l)$ be as described. Then for
$i_{t_s}\in I_l$ but $i_{t_s}\notin U\cup \{i_v\}$ we have
\begin{align}
C(I_l)&+L^*(U\mid I_l)-C(I_l\setminus \{i_{t_s}\})-L^*(U\mid I_l\setminus
\{i_{t_s}\}) \nonumber \\
&=-\sum_{k=v}^{s-1}\chi(i_{t_k}>j_{t_s}>i_v)a_{i_{t_k}}+\sum_{k=s+1}^{m-d}\chi(\overline{i_{t_k}>j_{t_s}>i_v})a_{i_{t_k}},
\end{align}
where
$\chi(\overline{i_{t_k}>j_{t_s}>i_v}):=1-\chi(i_{t_k}>j_{t_s}>i_v)$.
\end{lem}
We denote $-\sum_{k=v}^{s-1}\chi(i_{t_k}>j_{t_s}>i_v)a_{i_{t_k}}+\sum_{k=s+1}^{m-d}\chi(\overline{i_{t_k}>j_{t_s}>i_v})a_{i_{t_k}}$
by $g(i_{t_s})$.

\begin{lem}\label{lem-f1}
For $n\ge 2$, every term in the expansion of
$\prod_{s=1}^n\sum_{k\neq s}a(s,k)$ has $a(k,r)a(s,l)$ as a factor
for some $k,r,s,l$ satisfying $1\leqslant r\leqslant s<k\leqslant l\leqslant n$.
\end{lem}
\begin{proof}
Construct a matrix $A$ with $0$'s in the main diagonal as follows.
\begin{align*}
A=\begin{pmatrix}
0 & a(1,2) &\cdots &a(1,n) \\
a(2,1) &0  &\cdots &a(2,n) \\
\vdots &\vdots &\vdots &\vdots \\
a(n,1) &a(n,2) &\cdots &0
\end{pmatrix}.
\end{align*}
Then each term in the expansion of $\prod_{s=1}^n\sum_{k\neq s}a(s,k)$
corresponds to picking out one entry except for the $0$'s from each row of $A$. We prove by contradiction.

Suppose we choose $a(1,k_1)\ (k_1\geqslant 2)$ from the first row. Then
we can not choose $a(2,1)$, for otherwise $a(2,1)a(1,k_1)$ forms the
desired factor. Now from the second row, we have to choose
$a(2,k_2)\ (k_2\geqslant 3)$. It then follows that $a(3,1)$ and $a(3,2)$
can not be chosen, for otherwise $a(3,e)a(2,k_2),e=1,2$ forms the
desired factor. Repeat this discussion until the $n-1$st row, where
we have to choose $a(n-1,n)$. But then our $n$th row element
$a(n,e)$ (with $1\leqslant e\leqslant n-1$) together with $a(n-1,n)$ forms the
desired factor, a contradiction.
%
\end{proof}

The following factorization and cancelation lemma plays an important
role and it is our main discovery in this paper.
\begin{lem}\label{lem-factor}
For fixed set $U\neq I$ and integer $i_v\leqslant \min U$ we have the following
factorization {\small
\begin{align}\label{e-factor2} \sum_{I_l}(-1)^{l+d}q^{C(I_l)+L^*(U\mid I_l)} =(-1)^{\chi(\min
U \neq i_v)} q^{C(U\cup \{i_v\})+L^*(U\mid U\cup
\{i_v\})}\prod_{i_{t_s}\in I\setminus U\setminus
\{i_1,\ldots,i_v\}}\big(1-q^{g(i_{t_s})}\big),
\end{align}}
where $I_l$ ranges over all supersets of $U$ with the restriction
$\min I_l=i_v$.
Furthermore, if
there is no $s,t,u$ such that $1\leqslant s<t<u\leqslant m$ and
$j_t<i_s<j_u<i_t$, then
\begin{align}\label{e-cancel}
\prod_{i_{t_s}\in I\setminus U\setminus \{i_1,\ldots,i_v\}}\big(1-q^{g(i_{t_s})}\big)=0,
\end{align}
with the only exceptional case when $I\setminus U\setminus \{i_1,\ldots,i_v\}=\varnothing$.
\end{lem}

\begin{proof}
We prove this lemma in two parts.

\textbf{1.} Proof of \eqref{e-factor2}.

Notice that $I_l=U\cup \{i_v\}$ is the smallest set which satisfies
$\min I_l=i_v$ and $U\subseteq I_l$. So first we extract the common
factor $q^{C(U\cup \{i_v\})+L^*(U\mid U\cup \{i_v\})}$ from the
summation of \eqref{e-factor2}. Thus we need to calculate
$$C(I_l)+L^*(U\mid I_l)-C(U\cup \{i_v\})-L^*(U\mid U\cup \{i_v\}).$$
By Lemma \ref{lem-subtract} we have
\begin{align}\label{e-3}
C(I_l)&+L^*(U\mid I_l)-C(I_l\setminus \{i_{t_s}\})-L^*(U\mid I_l\setminus
\{i_{t_s}\})=g(i_{t_s}),
\end{align}
where $i_{t_s}\in I_l$ but $i_{t_s}\notin U\cup \{i_v\}$.
Thus iterating \eqref{e-3} we get
\begin{align}\label{e-factor3}
C(I_l)+L^*(U\mid I_l)-C(U\cup \{i_v\})-L^*(U\mid U\cup \{i_v\})
=\sum_{i_{t_s}\in I_l\setminus U\setminus \{i_v\}}g(i_{t_s}).
\end{align}
So extracting the common factor $q^{C(U\cup \{i_v\})+L^*(U\mid U\cup
\{i_v\})}$ from the left-hand side of \eqref{e-factor2} and by
\eqref{e-factor3} we have {\small \begin{align}\label{e-factor7}
\sum_{I_l}(-1)^{l+d}q^{C(I_l)+L^*(U\mid I_l)} = &q^{C(U\cup
\{i_v\})+L^*(U\mid U\cup \{i_v\})}\sum_{I_l}
(-1)^{l+d}q^{\sum_{i_{t_s}\in I_l\setminus U\setminus
\{i_v\}}g(i_{t_s})},
\end{align}}
where $I_l$ ranges over all supersets of $U$ with the restriction
$\min I_l=i_v$.

Second we prove the following factorization.
\begin{align}\label{e-factor6}
\sum_{I_l} (-1)^{l+d}q^{\sum_{i_{t_s}\in I_l\setminus U\setminus
\{i_v\}}g(i_{t_s})} =(-1)^{\chi(\min U \neq i_v)} \prod_{i_{t_s}\in
I\setminus U\setminus \{i_1,\ldots,i_v\}}\big(1-q^{g(i_{t_s})}\big),
\end{align}
where $I_l$ ranges over all supersets of $U$ and we restrict $\min
I_l=i_v$.

If $ \min U=i_v$, then the sign in the right-hand side of
\eqref{e-factor6} is positive. Every term in the expansion of the
right-hand side of \eqref{e-factor6} is of the form
$(-1)^{|G|}\prod_{i_{t_s}\in
G}q^{g(i_{t_s})}=(-1)^{|G|}q^{\sum_{i_{t_s}\in G}g(i_{t_s})}$,
where $G$ is a subset of $I\setminus U\setminus \{i_1,\ldots,i_v\}$.
Thus expanding the product of \eqref{e-factor6} we get
\begin{align}\label{e-change}
\prod_{i_{t_s}\in I\setminus U\setminus
\{i_1,\ldots,i_v\}}\big(1-q^{g(i_{t_s})}\big)=\sum_{G\subseteq
I\setminus U\setminus
\{i_1,\ldots,i_v\}}(-1)^{|G|}q^{\sum_{i_{t_s}\in G}g(i_{t_s})}.
\end{align}
Notice that
$I_l\setminus U\setminus \{i_v\}$ reduces to
$I_l\setminus U$ when $ \min U=i_v$.
Substitute $I_l\setminus U$ by $G'$ in the
left-hand side of \eqref{e-factor6}. Then $G'$ ranges over all
subsets of $I\setminus U\setminus \{i_1,\ldots,i_v\}$ if $I_l$
ranges over all supersets of $U$ with the restriction $\min I_l=i_v$.
Notice that $(-1)^{|G'|}=(-1)^{l-d}=(-1)^{l+d}$, thus the
left-hand side of \eqref{e-factor6} can also be written as the right
hand side of \eqref{e-change}. Hence \eqref{e-factor6} holds when
$\min U=i_v$. The case $\min U\neq i_v$ is similar.

Therefore \eqref{e-factor2} follows from \eqref{e-factor7} and \eqref{e-factor6}.

\textbf{2.} Under the assumption that
there is no $s,t,u$ such that $1\leqslant s<t<u\leqslant m$ and
$j_t<i_s<j_u<i_t$ we need to prove \eqref{e-cancel}.

If $\min I_l=\min U=i_v$, recall that $I\setminus
U=\{i_{t_1},\ldots,i_{t_{m-d}}\}$ and $t_1<\cdots<t_{m-d}$, then
$t_k=k$ for $k=1,\ldots,v-1$ and $t_v>v$. Thus $t_v\in I\setminus
U\setminus \{i_1,\ldots,i_v\}$. It follows that $\prod_{i_{t_s}\in
I\setminus U\setminus \{i_1,\ldots,i_v\}}\big(1-q^{g(i_{t_s})}\big)=
\prod_{s=v}^{m-d}\big(1-q^{g(i_{t_s})}\big)$.

If $\min I_l\neq \min U$, then ${t_v}=v$. It follows that $t_v\notin
I\setminus U\setminus \{i_1,\ldots,i_v\}$. Thus we have
$\prod_{i_{t_s}\in I\setminus U\setminus
\{i_1,\ldots,i_v\}}\big(1-q^{g(i_{t_s})}\big)=
\prod_{s=v+1}^{m-d}\big(1-q^{g(i_{t_s})}\big)$ and
$\chi(i_{t_v}>j_{t_s}>i_v)=\chi(i_{v}>j_{t_s}>i_v)=0$. In this case
$g(i_{t_s})$ reduces to
$$g(i_{t_s})=-\sum_{k=v+1}^{s-1}\chi(i_{t_k}>j_{t_s}>i_v)a_{i_{t_k}}+\sum_{k=s+1}^{m-d}\chi(\overline{i_{t_k}>j_{t_s}>i_v})a_{i_{t_k}}.$$
We only prove \eqref{e-cancel} when $\min I_l=\min U$, the case $\min I_l\neq \min U$ is similar.

We can write the left-hand side of \eqref{e-cancel} as $\prod_{s=v}^{m-d}\big(1-q^{g(i_{t_s})}\big)$ when $\min I_l=\min U$.
To prove
$\prod_{s=v}^{m-d}\big(1-q^{g(i_{t_s})}\big)=0$, it is sufficient to
prove $\prod_{s=v}^{m-d}g(i_{t_s})=0$.


Taking $a(s,k)=-\chi(i_{t_k}>j_{t_s}>i_v)a_{i_{t_k}}$ for $s>k$ and
$a(s,k)=\chi(\overline{i_{t_k}>j_{t_s}>i_v})a_{i_{t_k}}$ for $s<k$,
by the definition of $g(i_{t_s})$ we can write
$\prod_{s=v}^{m-d}g(i_{t_s})$ as $\prod_{s=v}^{m-d}\sum_{k\neq
s}a(s,k)$. By Lemma \ref{lem-f1} each term in the expansion of
$\prod_{s=v}^{m-d}g(i_{t_s})$ has a factor of the form
$-\chi(i_{t_r}>j_{t_k}>i_v)\chi(\overline{i_{t_l}>j_{t_s}>i_v})a_{i_{t_r}}a_{i_{t_l}}$,
where $v\leqslant r\leqslant s<k\leqslant l\leqslant m-d$. Thus
\begin{align}\label{e-factor5}
\prod_{s=v}^{m-d}g(i_{t_s})=\sum_{v\leqslant r\leqslant s<k\leqslant l\leqslant m-d}-\chi(i_{t_r}>j_{t_k}>i_v)\chi(\overline{i_{t_l}>j_{t_s}>i_v})a_{i_{t_r}}a_{i_{t_l}}\cdot \Delta,
\end{align}
where $\Delta$ is the product of some $a(s,k)$'s.

Next we prove each
$\chi(i_{t_r}>j_{t_k}>i_v)\chi(\overline{i_{t_l}>j_{t_s}>i_v})=0$ by contradiction
under the assumption that there is no $s,t,u$ such that $1\leqslant
s<t<u\leqslant m$ and $j_t<i_s<j_u<i_t$.

Suppose $\chi(i_{t_r}>j_{t_k}>i_v)\chi(\overline{i_{t_l}>j_{t_s}>i_v})=1$
for some $v\leqslant r\leqslant s<k\leqslant l\leqslant m-d$. Then
$\chi(i_{t_r}>j_{t_k}>i_v)=\chi(\overline{i_{t_l}>j_{t_s}>i_v})=1$.
By $\chi(i_{t_r}>j_{t_k}>i_v)=1$ we have
\begin{equation}\label{e-equation1}
i_{t_r}>j_{t_k}>i_v.
\end{equation}
By $\chi(\overline{i_{t_l}>j_{t_s}>i_v})=1$ we obtain
\begin{equation}\label{e-equation2}
i_{t_l}<j_{t_s}\quad \mbox{or} \quad j_{t_s}<i_v \quad \mbox{or}
\quad i_{t_l}<i_v.
\end{equation}
Since $l>v$, we have $t_l\geqslant l>v$ and $i_{t_l}>i_v$. Thus the last
inequality of \eqref{e-equation2} can not hold. Because $l>r$, $k>s$
and $i_{t_r}>j_{t_k}$ in \eqref{e-equation1}, we have
$i_{t_l}>i_{t_r}>j_{t_k}\geqslant j_{t_s}$. So the first inequality of
\eqref{e-equation2} can not hold too. Thus by \eqref{e-equation1}
and the middle inequality of \eqref{e-equation2} we obtain that if
$\chi(i_{t_r}>j_{t_k}>i_v)\chi(\overline{i_{t_l}>j_{t_s}>i_v})=1$
then $j_{t_s}<i_v<j_{t_k}<i_{t_r}$. It follows that
$j_{t_s}<i_v<j_{t_k}<i_{t_s}$ since $r\leqslant s$.
Because $v\leqslant s<k$, we have $v<t_v\leqslant t_s<t_k$. Thus for
$v<t_s<t_k$ the fact $j_{t_s}<i_v<j_{t_k}<i_{t_s}$ conflicts with
our assumption.
\end{proof}

\begin{lem}\label{lem-cancel}
If $U$ is of the form $\{i_{h},i_{h+1},\ldots,i_m\}$, then
\begin{align}\label{e-factor1}
q^{C(U)+L^*(U\mid U)}-q^{C(U\cup \{i_{h-1}\})+L^*(U\mid U\cup
\{i_{h-1}\})}=0.
\end{align}
\end{lem}
\begin{proof}
By the formula of $C(I_l)$ in \eqref{e-CIl} we have
\begin{align*}
C(U)+L^*(U\mid U)=1+a-\sum_{k\in U}a_k+\sum_{k=1}^{h-1}\big[N(i_{r_k},U)-N(i_{r_k},V_{k}^*)\big]a_{i_{r_k}},
\end{align*}
where $V_k^*=\{j_s>i_k\mid j_s\in V_1\cup \{j_{r_k}\}\}$ and $V_1=\{j_h,\ldots,j_{m}\}$ is the pairing set of $U$.
Since $U$ is of the form $\{i_{h},i_{h+1},\ldots,i_m\}$, we have
$i_{r_k}=i_k$ for $k=1,\ldots,h-1$. Hence $N(i_{r_k},U)=N(i_{r_k},V_{k}^*)=0$ for $k=1,\ldots,h-1$.
It follows that $C(U)+L^*(U\mid U)=1+a-\sum_{k\in U}a_k$.

Meanwhile
\begin{align*}
C&(U\cup \{i_{h-1}\})+L^*(U\mid U\cup\{i_{h-1}\})\nonumber \\
=&1+a-\sum_{k\in U}a_k-a_{i_{h-1}}+\sum_{k=1}^{h-2}\big[N(i_{r^{'}_k},U\cup\{i_{h-1}\})-N(i_{r^{'}_k},\overline{V_{k}^{*}})\big]a_{i_{r^{'}_k}}\nonumber \\
&-L^*(U\cup\{i_{h-1}\}\mid U\cup\{i_{h-1}\})+L^*(U\mid U\cup\{i_{h-1}\}),
\end{align*}
where $\overline{V_{k}^{*}}=\{j_s>i_k\mid j_s\in V_2\cup \{j_{r^{'}_k}\}\}$ and $V_2=\{j_{h-1},\ldots,j_{m}\}$.
Since $U\cup \{i_{h-1}\}$ is of the form $\{i_{h-1},i_{h},\ldots,i_m\}$, we have
$i_{r^{'}_k}=i_k$ for $k=1,\ldots,h-2$. Hence $N(i_{r^{'}_k},U\cup\{i_{h-1}\})=N(i_{r^{'}_k},\overline{V_{k}^{*}})=0$ for
$k=1,\ldots,h-2$. And by the definition of $L^{*}(T\mid
I)$ in \eqref{e-LT2} we have
$-L^*(U\cup\{i_{h-1}\}\mid U\cup\{i_{h-1}\})+L^*(U\mid U\cup\{i_{h-1}\})=a_{i_{h-1}}$.
Therefore $C(U\cup \{i_{h-1}\})+L^*(U\mid U\cup\{i_{h-1}\})$
has the same value as $C(U)+L^*(U\mid U)$.
\end{proof}

\subsection{Proof of the main theorem}

Having Lemma \ref{lem-factor} and Lemma \ref{lem-cancel}, we are ready to prove the main theorem.

\begin{proof}[Proof of Theorem \bf{\ref{t-m}}]
If $m=0$, then the theorem reduces to the $q$-Dyson Theorem. So we assume that $m\geqslant 1$.

Applying Theorem \ref{t-general} to the constant term in the left-hand side of \eqref{e-last} yields
{\small \begin{align}\label{e-expand}
\CT_{\mathbf{x}}&\left[\bigg(1+ \sum_{\varnothing \neq I_l\subseteq
I}(-1)^lq^{C(I_l)}\frac{x_{v_1}\cdots x_{v_l}}{x_{u_1}\cdots x_{u_l}}\bigg)
D_n(\mathbf{x},\mathbf{a},q)\right] \nonumber \\
&=\frac{(q)_{a}}{(q)_{a_0}\cdots(q)_{a_n}}
\bigg(1+\sum_{\varnothing\neq I_l\subseteq I}\sum_{\varnothing \neq
U\subseteq I_l}(-1)^{d+l}q^{C(I_l)+L^*(U\mid I_l)}
\frac{1-q^{\sum_{k\in U}a_k}}{1-q^{1+a-\sum_{k\in U}a_k}}\bigg),
\end{align}}
where $l=|I_l|$ and $d=|U|$.

Because $U$ is a subset of $I_l$, we have $\min I_l=i_v\leqslant \min U$.
By changing the summation order, the right-hand side of
\eqref{e-expand} can be rewritten as
\begin{equation}\label{e-new2}
\frac{(q)_{a}}{(q)_{a_0}\cdots(q)_{a_n}}\bigg(1+\sum_{\varnothing\neq
U\subseteq I}\sum_{i_{v}=i_1}^{\min U}\sum_{I_l}
(-1)^{d+l}q^{C(I_l)+L^*(U\mid I_l)} \frac{1-q^{\sum_{k\in
U}a_k}}{1-q^{1+a-\sum_{k\in U}a_k}}\bigg),
\end{equation}
where $I_l$ ranges over all supersets of $U$ with the restriction
$\min I_l=i_v$.

If $U\neq I$, then by Lemma \ref{lem-factor}, under the assumption
that there is no $s,t,u$ such that $1\leqslant s<t<u\leqslant m$ and
$j_t<i_s<j_u<i_t$ we have
\begin{align}\label{e-c2}
\sum_{I_l}(-1)^{l+d}q^{C(I_l)+L^*(U\mid I_l)}=0,
\end{align}
with the only exceptional case when $I\setminus U\setminus
\{i_1,\ldots,i_v\}=\varnothing$, where $I_l$ ranges over all
supersets of $U$ and we restrict $\min I_l=i_v$.

If $I\setminus U\setminus \{i_1,\ldots,i_v\}=\varnothing$, then $U$
is of the form $\{i_{h},i_{h+1},\ldots,i_{m}\}$ and $i_v$ is either
$i_h$ or $i_{h-1}$, and in this case $I_l=U$ or $I_l=U\cup \{i_{h-1}\}$
respectively. Thus by Lemma \ref{lem-cancel} we have
\begin{align}\label{e-c1}
q^{C(U)+L^*(U\mid U)}-q^{C(U\cup \{i_{h-1}\})+L^*(U\mid U\cup
\{i_{h-1}\})}=0.
\end{align}

By \eqref{e-c2} and \eqref{e-c1} the summands in \eqref{e-new2}
cancel with each other except for the summand when $U=I_l=I$. It
follows that \eqref{e-new2} reduces to
\begin{align}\label{e-2}
\frac{(q)_{a}}{(q)_{a_0}\cdots(q)_{a_n}}\Big(1+q^{C(I)+L^*(I\mid
I)}\frac{1-q^{\sum_{k\in I}a_k}}{1-q^{1+a-\sum_{k\in I}a_k}}\Big).
\end{align}
By the formula of $C(I_l)$ in \eqref{e-CIl} we get
$C(I)=1+a-\sum_{k\in I}a_k-L^*(I\mid I)$. Substituting $C(I)$ into
\eqref{e-2} and multiplying the equation by $1-q^{1+a-\sum_{k\in I}a_k}$ we can obtain the right-hand side of \eqref{e-last}.
\end{proof}

\section{Remark}

%
%

If there exist some $s,t,u$ such that $s<t<u$ and $j_t<i_s<j_u<i_t$,
then our main theorem does not lead to the desired cancelations. As
stated in Section 4.1, we can solve for $A_T$ such that the constant
term of $\sum_{T}A_T \frac{x_{v_1}\cdots x_{v_l}}{x_{u_1}\cdots
x_{u_l}} D_n(\mathbf{x},\mathbf{a},q)$ has the desired cancelations.
However, experiments show that there is no nice form for $A_T$ in this situation.

Another possibility to let the $u$'s and the $v$'s be not paired up.
Some of the cases can be established by applying the operator $\pi$
defined in Section 2 to our main theorem. But not all the un-paired
up cases can be obtained in this way.

\vspace{.2cm} \noindent{\bf Acknowledgments.}
Some of the results in this paper were obtained in the Center for Combinatorics of
Nankai University when I studied there.
I
would like to acknowledge the helpful guidance of my supervisor
William Y.C. Chen.
I am very grateful to Guoce Xin, for his guidance, suggestions and help.
I thank Lun Lv for helping me check the errors in my paper.

\section{Appendix: Proof of Theorem \ref{t-general} }

\begin{proof}
By the definition of $\pi$, it is easy to deduce that
\begin{equation}\label{e-def-pi}\pi^k x_i=\left\{ \begin{aligned}
         &x_{i+k},\quad &for& \quad i+k\leqslant n; \\
                  &x_{i+k-n-1}/q, \quad &for& \quad i+k>n.
                          \end{aligned} \right.\end{equation}
Iterating Lemma \ref{l-pi} $n-i_1+1$ times, i.e., acting with
$\pi^{n-i_1+1}$, we obtain

\begin{align}
&\CT_{\mathbf{x}}\frac{x_{j_1}\cdots
x_{j_{m}}}{x_{i_1}\cdots
x_{i_m}}D_n(\mathbf{x},\mathbf{a},q) \nonumber \\
=&\CT_{\mathbf{x}}\frac{\prod_{l=1}^tx_{j_l+n-i_1+1}\prod_{l=t+1}^{m}x_{j_l-i_1}q^{-(m-t)}}{x_0x_{i_2-i_1}\cdots
x_{i_m-i_1}q^{-m}} D_n(\mathbf{x},(b_0,\ldots,b_n),q),
\end{align}
where
\begin{equation}\label{e-def-b}
b_{k}=\left\{ \begin{aligned}
         &a_{k+i_1},\quad &for& \quad k=0,\ldots,n-i_1; \\
                  &a_{k-(n-i_1+1)}, \quad &for& \quad k=n-i_1+1,\ldots,n.
                          \end{aligned} \right. \end{equation}

To apply Theorem \ref{t-main-thm}, we define
$\widetilde{I}=\{0,i_2-i_1,\ldots,i_m-i_1\}$, and
$\widetilde{J}^-=\{j_1+n-i_1+1,\cdots,j_t+n-i_1+1\}$,
$\widetilde{J}^+=\{j_{t+1}-i_1,\ldots,j_m-i_1\}$,
$\widetilde{J}=\widetilde{J}^-\cup \widetilde{J}^+$. Then by Theorem
\ref{t-main-thm} we have
\begin{align*}
\CT_{\mathbf{x}}\frac{x_{j_1}\cdots x_{j_{m}}}{x_{i_1}\cdots
x_{i_m}}D_n(\mathbf{x},\mathbf{a},q) 
=&q^{t}\frac{(q)_{a}}{(q)_{a_0}\cdots(q)_{a_n}}\sum_{
\varnothing\neq \widetilde{T}\subseteq \widetilde{I}}(-1)^dq^{L(\widetilde{T}\mid \widetilde{I})}
\frac{1-q^{\sum_{k\in \widetilde{T}}b_k}}{1-q^{1+a-\sum_{k\in \widetilde{T}}b_k}},
\end{align*}
where $|\widetilde{T}|=d$ and
\begin{align}\label{e-LT'}
L(\widetilde{T}\mid \widetilde{I})=\sum_{k=0}^n\big[N(k,\widetilde{I})-N(k,\widetilde{J})\big]\widetilde{w}_k,
\end{align}
in which $\widetilde{w}_k$ is $b_k$ if $k\notin \widetilde{T}$ and
$0$ otherwise.

There is a natural one-to-one correspondence between $I$ and
$\widetilde{I}$: $I\mathop{\longrightarrow}\limits^{f}
\widetilde{I}$, $f(a)=a-i_1$, $a\in I$. This correspondence clearly
applies between their subsets $T$ and $\widetilde{T}$.

Since the largest element in $\widetilde{T}$ is not larger than
$i_m-i_1$ and $i_m-i_1\leqslant n-i_1$, by the definition of $b_k$ we
have
\begin{equation*}\sum_{k\in \widetilde{T}}b_k=\sum_{k\in \widetilde{T}}a_{k+i_1}=\sum_{k\in T}a_k.\end{equation*}

Next we have to rewrite \eqref{e-LT'} in terms of $w_k$, $N(k,I)$
and $N(k,J)$ to get $L^*(T\mid I)$.

Because the largest element in $\widetilde{I}$ is $i_{m}-i_1\leqslant
n-i_1$, so if $k>n-i_1$ then $k\notin \widetilde{T}$. It follows
that
\begin{equation}\label{e-trans1}
\widetilde{w}_k=b_k=a_{k-(n-i_1+1)}.
\end{equation}
If $k\leqslant n-i_1$, then
\begin{equation} \label{e-trans2}
\widetilde{w}_k=\left\{ \begin{aligned}
         &b_k=a_{k+i_1}, \quad &\mbox{if} \quad k\notin \widetilde{T}; \\
                  &0, \quad &\mbox{if} \quad k\in \widetilde{T},
                          \end{aligned} \right.
                          \end{equation}
which is in fact $w_{k+i_1}$.

It is straightforward to check that
\begin{equation}\label{e-trans3}
N(k,\widetilde{I})=N(k+i_1,I),
\end{equation}
\begin{equation}\label{e-trans4}
N(k,\widetilde{J}^-)=N(k-(n-i_1+1),J^-),\quad N(k,\widetilde{J}^+)=N(k+i_1,J^+),
\end{equation}
\begin{equation}\label{e-trans5}
N(k,\widetilde{J})=N(k,\widetilde{J}^-)+N(k,\widetilde{J}^+).
\end{equation}

Substituting \eqref{e-trans1} and \eqref{e-trans2} into  \eqref{e-LT'} we have
\begin{align*}
L(\widetilde{T}\mid \widetilde{I})=&\sum_{k=0}^{n-i_1}\big[N(k,\widetilde{I})-N(k,\widetilde{J})\big]w_{k+i_1}+\sum_{k=n-i_1+1}^{n}\big[N(k,\widetilde{I})-N(k,\widetilde{J})\big]a_{k-(n-i_1+1)}.
\end{align*}
By \eqref{e-trans3}--\eqref{e-trans5} the above equation becomes
{\small\begin{align}\label{e-trans6}
L(\widetilde{T}\mid \widetilde{I})=&\sum_{k=0}^{n-i_1}\big[N(k+i_1,I)-N(k-(n-i_1+1),J^-)-N(k+i_1,J^+)\big]w_{k+i_1} \nonumber \\
&+\sum_{k=n-i_1+1}^{n}\big[N(k+i_1,I)-N(k-(n-i_1+1),J^-)-N(k+i_1,J^+)\big]a_{k-(n-i_1+1)}.
\end{align}}

If $k\in [0,n-i_1]$ then $k-(n-i_1+1)<0$. Thus $N(k-(n-i_1+1),J^-)=0$.

If $k\in [n-i_1+1,n]$ then $k+i_1>n$. Thus $N(k+i_1,I)=m$ and $N(k+i_1,J^+)=m-t$.

Therefore \eqref{e-trans6} reduces to
{\small\begin{align*}\label{e-trans7}
&L(\widetilde{T}\mid \widetilde{I})\nonumber \\
&=\sum_{k=0}^{n-i_1}\big[N(k+i_1,I)-N(k+i_1,J^+)\big]w_{k+i_1}
+\sum_{k=n-i_1+1}^{n}\big[t-N(k-(n-i_1+1),J^-)\big]a_{k-(n-i_1+1)}\nonumber \\
&=\sum_{k=i_1}^{n}\big[N(k,I)-N(k,J^+)\big]w_{k}
+\sum_{k=0}^{i_1-1}\big[t-N(k,J^-)\big]a_{k}.
\end{align*}}
Then we obtain
\begin{align*}
L^*(T\mid I)=t+L(\widetilde{T}\mid \widetilde{I})=t+\sum_{k=i_1}^{n}\big[N(k,I)-N(k,J^+)\big]w_{k}
+\sum_{k=0}^{i_1-1}\big[t-N(k,J^-)\big]a_{k}.
\end{align*}
\end{proof}


\begin{thebibliography}{10}

\bibitem{andrews1975}
G. E. Andrews, \emph{Problems and prospects for basic hypergeometric
functions}, in \emph{Theory and Application of Special Functions},
ed. R. Askey, Academic Press, New York, 1975, pp. 191--224.


\bibitem{breg}
D.~M. Bressoud and I.~P. Goulden, \emph{Constant term identities
extending the
  {$q$}-{D}yson theorem}, Trans. Amer. Math. Soc. \textbf{291} (1985),
  203--228.

\bibitem{dyson1962}
F. J. Dyson, \emph{Statistical theory of the energy levels of complex systems I},
J. math. Phys. \textbf{3} (1962), 140--156.

\bibitem{gess-xin2006}
I. M. Gessel and G. Xin, \emph{A short proof of the
Zeilberger-Bressoud $q$-Dyson theorem}, Proc. Amer. Math. Soc.
\textbf{134} (2006), 2179--2187.

\bibitem{good1}
I.~J. Good, \emph{Short proof of a conjecture by {D}yson}, J. Math.
Phys.
  \textbf{11} (1970), 1884.

\bibitem{gunson}
J.~Gunson, \emph{Proof of a conjecture by {D}yson in the statistical
theory of
  energy levels}, J. Math. Phys. \textbf{3} (1962), 752--753.


\bibitem{kadell1998}
K. W. J. Kadell, \emph{Aomoto's machine and the Dyson constant term
identity}, Methods Appl. Anal. \textbf{5} (1998), 335--350.

\bibitem{Xin}
L. Lv, G. Xin and Y. Zhou, \emph{A Family of $q$-Dyson style
constant term identities}, J. Combin. Theory. Ser. A \textbf{116}
(2009), 12--29.

\bibitem{Xin2}
L. Lv, G. Xin and Y. Zhou, \emph{Two coefficients of the Dyson product},
Electro. J. Combin. \textbf{15} (2008), R36, 11 pp.

\bibitem{macdonald}
I.~G.~Macdonald, \emph{Some conjectures for root systems}, SIAM J.
Math. Anal. \textbf{13} (1982), 988--1007.

\bibitem{sills2006}
A. V. Sills, \emph{Disturbing the Dyson conjecture, in a generally
GOOD way}, J. Combin. Theory Ser. A \textbf{113} (2006), 1368--1380.

\bibitem{sills-zeilberger}
A. V. Sills and D. Zeilberger, \emph{Disturbing the Dyson conjecture
(in a Good way)}, Experiment. Math. \textbf{15} (2006), 187--191.



\bibitem{stembridge1987}
J.~R. Stembridge, \emph{First layer formulas for characters of
$SL(n,\mathbb{C})$}, Trans. Amer. Math. Soc.  \textbf{299} (1987),
319--350.

\bibitem{stembridge-qdyson}
J.~R. Stembridge, \emph{A short proof of {M}acdonald's conjecture
for the root
  systems of type {$A$}}, Proc. Amer. Math. Soc. \textbf{102} (1988),
  777--786.



\bibitem{wilson}
K.~G. Wilson, \emph{Proof of a conjecture by {D}yson}, J. Math.
Phys. \textbf{3} (1962), 1040--1043.

\bibitem{zeil1982}
D.~Zeilberger, \emph{A combinatorial proof of Dyson's conjecture},
Discrete Math. \textbf{41} (1982), 317--321.

\bibitem{zeil-bres1985}
D. Zeilberger and D. M. Bressoud, \emph{A proof of Andrews'
$q$-Dyson conjecture}, Discrete Math. \textbf{54} (1985),
201--224.

\end{thebibliography}
\end{document}